\documentclass{amsart}
\font\goth=eusm10

\newcommand\Ii{\hbox{\goth I}}

\newcommand\ZZ{\mathbb{Z}}
\newcommand\Q{\mathbb{Q}}
\newcommand\C{\mathbb{C}}
\newcommand\R{\mathbb{R}}

\newcommand\Oc{\hbox{\goth O}}

\newcommand\Pq{\mathbb{P}^4}

\newcommand\Pt{\mathbb{P}^3}
\newcommand\Pd{\mathbb{P}^2}
\newtheorem{theorem}{Theorem}

\newtheorem{proposition}{Proposition}
\newtheorem{definition}{Definition}

\newtheorem{remark}{Remark}

\begin{document}

\title{Some results of regularity for Severi varieties of projective 
surfaces}

\author{F. Flamini}

\date{preprint submitted on March 22, 1999}


\curraddr{Dipartimento di Matematica, Terza Universita' degli Studi di Roma 
- \\''Roma Tre''-  L. S. L. Murialdo, 1 - \\ 00146 Roma - Italy}
\email{flamini@matrm3.mat.uniroma3.it}

\begin{abstract}
For a linear system $\mid C \mid$ on a smooth projective surface 
$S$, whose general member is a smooth, irreducible
curve, the \emph{Severi variety} $V_{\mid C \mid, \delta}$ is the
locally closed subscheme of $\mid C \mid$ which parametrizes irreducible curves 
with only $\delta$
nodes as singularities. In this paper we give numerical conditions
on the class of divisors and upper bounds on $\delta$, ensuring
that the correspondent Severi variety is smooth of codimension $\delta$. 
Our result generalizes what is proven in \cite{C-S} and \cite{G-L-S II}.  
We also consider examples of smooth Severi varieties on surfaces of general 
type in $\Pt$ which contain a line.
\end{abstract}

\thanks{The author is a member of GNSAGA-CNR}
\thanks{This work is part of my Ph.D. thesis. I wish to 
thank many people who helped me during its preparation, especially 
L. Chiantini, D. Franco, M. Franciosi, A. F. Lopez, G. Pareschi, E. Sernesi 
and A. Verra. My special thanks to my Ph.D. advisor E. Sernesi for having 
introduced me in such an interesting research area.}
 
\maketitle

\section*{Introduction.}    
     
Nodal curves play a central role in the subject of singular curves. The 
definition of the \emph{Severi variety} of irreducible, nodal curves on any 
smooth, projective surface is standard. For a given effective divisor 
$C \in Div(S)$, let $\mid C \mid$ denote the linear system associated to the 
line bundle $\Oc_S(C)\in Pic(S)$. If we suppose that the generic element 
of $\mid C \mid$ is a smooth, irreducible curve, it makes sense to consider 
the subscheme $ V_{\mid C \mid, \; \delta}$ of $\mid C \mid$, which 
parametrizes all curves $ C' \in \mid C \mid$ that are irreducible and have 
only $\delta$ nodes as singular points. It is well known that such a 
subscheme is locally closed in the projective space $\mid C \mid$.     
     
In \cite{Sev}, Anhang F, Severi studied some properties of the variety $V_{d, \; g}$, defined as the closure of the locus consisting of reduced and 
irreducible plane curves, of geometric genus $g$, in the projective space parametrizing plane curves of degree $d$. $V_{d, \; g}$ contains, as an open dense subscheme, the locus $V_{d, \; \delta}$ corresponding to irreducible curves having only $\delta$ nodes as singularities. $V_{d, \; g}$ is classically known as the Severi variety of plane curves of given genus and degree. He proved that, for every $ d \geq 3$ and 
$ 0\leq \delta \leq \left( \begin{array}{c}
                                d-1 \\ 2
                         \end{array} \right)$, $ V_{\mid dL \mid, 
\; \delta}$ is non-empty and everywhere smooth of codimension $\delta$ 
in $\mid dL \mid$, where $L$ denotes a line in $\Pd$. Only after more than 
60 years, Harris completed, in \cite{H}, Severi's proof of the 
irreducibility of the Severi variety $V_{d, \; g}$ of the projective plane, 
by showing that the dense open subset $V_{d, \; \delta}$ is connected. 

With abuse of terminology, we shall use in the sequel the word {\it Severi variety} for the locally closed subscheme $V_{\mid C \mid, \; \delta}$ of the linear series $\mid C \mid$ on a projective surface $S$. 

In recent times, there have been many results on this subject and in many 
directions. In fact, one may study several problems concerning Severi 
varieties.

Existence problems are covered by recent investigations in the case of 
Del Pezzo surfaces \cite{G-L-S I} or K3 surfaces \cite{Ch}; 
on the other hand, we 
have a complete computation of the degree of Severi varieties in the plane, 
treated in \cite{C-H}. There are some results even on the irreducibility 
problem,
 contained in \cite{C-C}, 
for general surfaces of $\Pt$ of degree $d\geq 8$. 

A natural approach to the dimension problem is to use deformation theory of 
nodal curves. Apart from the Severi classical result, whose proof can be 
extended to some rational or ruled surfaces and to K3 surfaces, 
information about regularity of $ V_{\mid C \mid, \; \delta}$ on a surface 
$S$ of general type can be obtained by studying suitable rank 2 vector 
bundles on $S$. The first who used such approach were Chiantini and Sernesi. 
In \cite{C-S} they found an upper bound on $\delta$, ensuring that the 
family 
$ V_{\mid C \mid, \; \delta}$ is smooth of codimension $\delta$ in the 
projective space $\mid C \mid$. Their proof focused on surfaces such that 
$\mid K_S \mid$ is ample and $C$ is a divisor which is numerically 
equivalent to $pK_S$, $p \in \Q^+$ and $p \geq 2$.

An improvement of this result is given in \cite{G-L-S II}. The authors 
generalized this approach in two directions. In fact, they allowed 
arbitrary singularities and they weakened the assumption of $K_S$ being 
ample, so that $S=\Pd$ is included. Their assumptions are: $C$, $C-K_S$ 
ample divisors and $C^2 \geq K_S^2$; moreover, some numerical hypotheses 
are made, which imply, in the case of nodes, that $(C-2K_S)^2 >0$ and 
$C(C-2K_S)>0$.

In this paper, we give a purely numerical criterion to prove the regularity 
of $ V_{\mid C \mid, \; \delta}$ , provided that $\delta$ is less than a 
suitable upper bound. More precisely, we prove the following

$ \!\!$
\linebreak
{\bf Theorem} \emph{Let $S$ be a smooth, projective surface and $C$ be a 
smooth, irreducible divisor on $S$. Suppose that:
\begin{enumerate}
\item $(C-2K_S)^2 >0$ and $C(C-2K_S)>0$;
\item either\begin{displaymath}\begin{array}{lclc}
K_S^2>-4 & if & C(C-2K_S)\geq 8, & or \\
K_S^2\geq 0 & if & 0 < C(C-2K_S)<8. &  
\end{array}
\end{displaymath}
\item $CK_S \geq 0$;
\item $H(C,K_S)< 4(C(C-2K_S)-4)$, where $H(C,K_S)$ is the $Hodge \; 
number$ (see def. \ref{def:2}) of $C$ and $S$;
\item
\begin{displaymath}\begin{array}{lclc}
\delta \leq \frac{C(C-2K_S)}{4}-1& if & C(C-2K_S)\geq 8, & or \\
\delta < \frac{C(C-2K_S)+\sqrt{C^2(C-2K_S)^2}}{8}& if & 0 < C(C-2K_S)<8. &    
\end{array}
\end{displaymath}
\end{enumerate}
Then, if $C' \in \mid C\mid$ is a reduced, irreducible curve 
with only $\delta$ nodes as singular points and if $N$ denotes the 
$0$-dimensional scheme of nodes in $C'$, in the above hypotheses $N$ 
imposes independent conditions to $\mid C \mid$, i.e. the Severi variety 
$ V_{\mid C \mid, \; \delta}$ is smooth of codimension $\delta$ at $C'$}.

We shall give some examples which show that such a result really 
generalizes the ones recalled before. Moreover, we can obtain some results 
on the Severi varieties of surfaces in $\Pt$ which are elements of 
a component of the Noether-Lefschetz locus, consisting of surfaces which 
contain a line. This is related to some results contained in 
\cite{C-L}, where the question of algebraic hyperbolicity 
for surfaces $S$ in $\Pt$ and in $\Pq$ is treated.

In the sequel, we shall work in the category of $\C-schemes$. We will 
denote by $\equiv$ the linear equivalence of divisors, whereas 
$\equiv_{num}$
 shall denote the numerical equivalence of divisors.

\section*{Preliminaries}

Let $S$ be a projective, non-singular algebraic surface and $\mid D \mid$ a 
complete linear system on $S$ whose general member is a smooth, irreducible 
curve. If $p_a(D)$ denotes the $arithmetic \; genus$ of $D$ then, by the 
adjunction formula,
$$p_a(D) = \frac{D(D+K_S)}{2}+1.$$For a given $\delta \geq 1$, suppose that 
$ V_{\mid D \mid, \; \delta}$ is non-empty. Let $C \in  V_{\mid C \mid, \; 
\delta}$ and let $N$ be the scheme consisting of the $\delta$ nodes of $C$. 
The $ geometric \;genus$ of $C$ is $g=p_g(C)=p_a(C)-\delta.$

We know that the Zariski tangent space of $\mid D \mid$ at the point 
$[C]$, parametrizing $C$, is isomorphic 
to $$H^0(S, \Oc_S(D))/< C >,$$ whereas the Zariski tangent space to 
$ V_{\mid D \mid, \; \delta}$ at $[C]$ is

$$ T_C(V_{\mid D \mid, \; \delta}) \cong H^0(S, \Ii_N(D))/< C >,$$ where 
$\Ii_N \subset \Oc_S$ denotes the ideal sheaf of the subscheme $N$ of $S$ 
(see, for example, \cite{S}). The relative obstruction space is a subspace 
of 
$H^1(S, \Ii_N(D))$. In particular, $N$ imposes independent conditions to 
$\mid D \mid$ if and only if

$$dim(V_{\mid D \mid, \; \delta}) = dim \; T_C(V_{\mid D \mid, \; \delta})
= dim \mid D \mid - \delta$$ at $C$. In this case, $V_{\mid D \mid, \; 
\delta}$ is smooth of the $expected \; dimension$ at [$C$]. 
The component containing [$C$] will be called 
\emph{regular} at [$C$]. Otherwise, it is said to be a \emph{superabundant 
component}. The regularity property is very strong, since it implies that 
the nodes of $C$ can be independently smoothed (see \cite{C-C} or \cite{C-S}).

\begin{definition}\label{def:1}

\normalfont{Let $S$ be a smooth, projective surface and $C \in Div(S)$. 
$C$ is said to be a \emph{nef divisor} if $CF\geq 0$, for each effective 
divisor $F$. A nef divisor $C$ is called $big$ if $C^2>0$.}          
\end{definition}

We recall that, by the \emph{Kleiman criterion} (see 
\cite{Hart}), $C$ is nef if and only if it is in 
the closure of the ample divisor cone of $S$.

\begin{definition}\label{def:2}
\normalfont{Let $S$ be a smooth, projective surface and $C \in Div(S)$. 
We shall denote by $H(C, K_S)$ the \emph{Hodge number} of $C$ and $S$, 
defined by $$ H(C,K_S):= (CK_S)^2 - C^2K_S^2.$$}
\end{definition}
The \emph{Index theorem} (see \cite{B-P-V}, page 120) ensures us that this 
number is a non-negative real number.

\section*{The main result}

In the previous section we have recalled all definitions and properties 
needed to prove our principal result. We are now able to state the following

\begin{theorem}\label{main} Let $S$ be a smooth, projective surface and 
$C$ be a smooth, irreducible divisor on $S$. Suppose that:
\begin{enumerate}
\item $(C-2K_S)^2 >0$ and $C(C-2K_S)>0$;
\item\begin{displaymath}\begin{array}{llclc}
(i) & K_S^2>-4 & if & C(C-2K_S)\geq 8, & or \\
(ii) & K_S^2\geq 0 & if & 0 < C(C-2K_S)<8. &  
\end{array}
\end{displaymath}
\item $CK_S \geq 0$;
\item $H(C,K_S)< 4(C(C-2K_S)-4)$, where $H(C,K_S)$ is the $Hodge \; number$ 
of $C$ and $S$ (def.\ref{def:2});
\item
\begin{displaymath}\begin{array}{llclc}
(i) & \delta \leq \frac{C(C-2K_S)}{4}-1& if & C(C-2K_S)\geq 8, & or \\
(ii) & \delta < \frac{C(C-2K_S)+\sqrt{C^2(C-2K_S)^2}}{8}& if & 0 < C(C-2K_S)
<8. &    
\end{array}
\end{displaymath}
\end{enumerate}
Then, if $C' \in \mid C\mid$ is a reduced, irreducible curve with only 
$\delta$ nodes as singular points and if $N$ denotes the $0$-dimensional 
scheme of nodes in $C'$, in the above hypotheses $N$ imposes independent 
conditions to $\mid C \mid$, i.e. the Severi variety $ V_{\mid C \mid, \; 
\delta}$ is smooth of codimension $\delta$ at $C'$.
\end{theorem}
\begin{proof}
For the sake of simplicity we will write $K$, instead of $K_S$, to denote a 
canonical divisor of $S$. By contradiction, assume that $N$ does not impose 
independent conditions to $\mid C \mid$. Let $N_0 \subset N$ be a minimal 
0-dimensional subscheme of $N$ for which this property holds and let 
$\delta_0= \mid N_0 \mid$. This means that $H^1(S, \; \Ii_{N_0}(C)) \neq 0$ 
and that $N_0$ satisfies the \emph{Cayley-Bacharach} condition 
(see, for example,
 \cite{G-H}). Therefore, a non-zero element of $H^1(\Ii_{N_0}(C))$ 
corresponds 
to a non-trivial rank 2 vector bundle $E \in Ext^1(\Ii_{N_0}(C-K), \Oc_S)$; 
so, one can consider the following exact sequence
\begin{equation}
0 \to \Oc_S\to E\to \Ii_{N_0}(C-K)\to 0. \label{eq:1.1}
\end{equation}This implies that

$$c_1(E)=C-K,\; c_2(E)= \delta_0 \leq \delta,$$i.e. 
\begin{equation}
 c_1(E)^2 - 4 c_2(E)=(C-K)^2 - 4 \delta_0. \label{eq:1.2}
\end{equation}We now want to compute $(\ref{eq:1.2})$ in cases $5.(i)$ 
and $5.(ii)$.

In the first one, $$(C-K)^2 - 4 \delta_0 \geq (C-K)^2 - 4 \delta= C^2 - 
2CK -4 + 4 + K^2 - 4 \delta \geq  K^2 + 4 > 0,$$ by $2(i)$.

In the other case, using $5.(ii)$ and the Index Theorem, $$(C-K)^2 - 4 
\delta_0 \geq (C-K)^2 - 4 \delta= C^2 - 2CK + K^2 - 4 \delta >  K^2 \geq 
0,$$ since we supposed $2(ii)$.

In both cases, the vector bundle $E$ is $Bogomolov-unstable$ (see
\cite{Bog} or \cite{Reid}), i.e. there 
exist $M, \; B \in Div(S)$ and a 0-dimensional scheme $Z$ (possibly empty) 
such that
\begin{equation}
0 \to \Oc_S(M)\to E\to \Ii_{Z}(B)\to 0 \label{eq:1.3}
\end{equation}holds and $(M-B) \in N(S)^+$. We recall that $N(S)^+$ 
denotes the ample divisor cone of $S$. This means that 
\begin{eqnarray}
(M-B)^2 & > 0,& \nonumber \\
     &   &    \\       \label{eqnarray:1.4}
(M-B)H & >0, & \forall \;H \; ample \; divisor. \nonumber
\end{eqnarray}The exact sequence $(\ref{eq:1.3})$ ensures us that $H^0 
(E(-M)) \neq 0$. If we consider the tensor product of the exact sequence 
$(\ref{eq:1.1})$ by $\Oc_S(-M)$, we get
\begin{equation}
0 \to \Oc_S(-M)\to E(-M)\to \Ii_{N_0}(C-K-M)\to 0. \label{eq:1.5}
\end{equation}We state that $H^0 (\Oc_S(-M))=0$; otherwise, $-M$ would 
be an effective divisor, therefore $-MH>0$ for each ample divisor $H$. 
From $(\ref{eq:1.3})$, it follows that $c_1(E)=M+B$, so, by $(\ref{eq:1.1})$
 and $(4)$, 
\begin{equation}\label{eq:1.6}
M-B = 2M-C+K \in N(S)^+.
\end{equation} Thus, for every ample divisor $H$,

\begin{equation}
MH > \frac{(C-K)H}{2}. \label{eq:1.7}
\end{equation}Furthermore, from hypotheses $1.$ and $3.$, it immediately 
follows that $C(C-K)>0$ and $C^2>0$. Indeed, $C(C-K)=C(C-2K)+CK>0$ and 
$C^2>CK \geq 0$. Since $C$ is irreducible, this implies that $C$ is a nef 
divisor; from $(\ref{eq:1.7})$ and from Kleiman's criterion, we get
\begin{equation}
MC \geq \frac{(C-K)C}{2}.\label{eq:1.8}
\end{equation}It follows that $-MC<0$ so, since $C$ is nef, $-M$ can not 
be effective.

If we consider the cohomological exact sequence associated to 
$(\ref{eq:1.5})$, we deduce that there exists a divisor $\Delta \in \mid 
C-K-M \mid$ s.t. $N_0 \subset \Delta$. If the irreducible nodal curve $C' 
\in \mid C \mid$, whose sets of nodes is $N$, were component of $\Delta$, 
then $-M-K$ would be an effective divisor. By applying $(\ref{eq:1.8})$ and 
by using the fact that $C(C-K)>0$ and hypothesis $3.$, one determines
\vspace{-10pt}
$$C'(-M-K)= C(-M-K)= -CK-CM \leq - CK -\frac{(C-K)C}{2}=$$
$$=- \frac{(C+K)C}{2}=-\frac{KC}{2} - \frac{C^2}{2}< -CK \leq 0,$$
which contradicts the effectiveness of $-M-K$, since $C$ is nef.

Bezout's theorem implies that
\begin{equation}
C' \Delta =C'(C-K-M)\geq 2 \delta_0. \label{eq:1.9}
\end{equation}On the other hand, taking $M$ maximal, we may further assume 
that the general section of $E(-M)$ vanishes in a $2$-codimensional locus 
$Z$ of $S$. Thus, $c_2(E(-M))= deg(Z) \geq 0$. By standard computations on 
Chern classes, we obtain

\vspace{-10pt}
$$c_2(E(-M))= c_2(E) + M^2+ c_1(E)(-M) = \delta_0 + M^2 - M(C-K),$$which 
implies 
\begin{equation}
\delta_0 \geq M(C-K-M). \label{eq:1.10}
\end{equation}

By applying the Index theorem to the divisor pair $(C,\; 2M-C+K)$, we get
\begin{equation}
C^2(2M-C+K)^2 \leq (C(C-K)-2C(C-K-M))^2. \label{eq:1.11}
\end{equation}

From $(\ref{eq:1.9})$ and from the fact that $C(C-K)$ is positive, it 
follows that 
\begin{equation}
C(C-K)-2C(C-K-M) \leq C(C-K) - 4 \delta_0.\label{eq:1.12}
\end{equation} We observe that the left side member of (\ref{eq:1.12}) is 
non-negative, since $C(C-K)-2C(C-K-M) = C(2M-C+K)$, where $C$ is effective 
and, by ($\ref{eq:1.6}$), $2M-C+K \in N(S)^+$. Thus, ($\ref{eq:1.12}$) 
still holds when we square both its members and, together with 
($\ref{eq:1.11}$), this gives
\begin{equation}
C^2(2M-C+K)^2 \leq (C(C-K)-4 \delta_0)^2.\label{eq:1.13}
\end{equation}

On the other hand, using $(\ref{eq:1.10})$, we get
$$(2M-C+K)^2= 4(M- \frac{(C-K)}{2})^2 =$$ $$ =(C-K)^2-4(C-K-M)M \geq 
(C-K)^2- 4 \delta_0,$$i.e.
\begin{equation}
(2M-C+K)^2 \geq (C-K)^2 - 4 \delta_0. \label{eq:1.14}
\end{equation}Putting together $(\ref{eq:1.13})$ and $(\ref{eq:1.14})$, 
we get

\begin{equation}
F(\delta_0):= 16 \delta_0^2 - 4 C(C-2K) \delta_0 + (CK)^2 - C^2 K^2 \geq 
0.\label{eq:1.15}
\end{equation} Summarizing, the assumption on $N$, stated at the beginning,
 implies $(\ref{eq:1.15})$\footnote{We remark that, in the case of nodes,
 this condition is the 
same of \cite{G-L-S II}; moreover, their hypotheses $(1.2)$ and $(1.3)$ 
coincide 
in the case of nodes and become $F(\delta_0) <0$.}. We want to show that our numerical hypotheses 
hold if and only if the opposite inequality is satisfied. To this aim, 
observe that the discriminant of the equation $F(\delta_0)=0$ is 
$16 C^2(C-2K)^2$, so, by hypotheses $1.$ and $3.$, it is positive. 
The inequality $F(\delta_0)<0$ is verified iff $\delta_0 \in 
( \alpha(C,K), \; \beta(C, \; K))$, where

$$ \alpha(C,K)=\frac{C(C-2K)-\sqrt{C^2(C-2K)^2}}{8} \in \R \; and$$ 
$$ \beta(C,K)=\frac{C(C-2K)+\sqrt{C^2(C-2K)^2}}{8}\in \R\; ;$$ we have
 to show that, with our numerical hypotheses, $\delta_0 \in 
( \alpha(C,K), \; \beta(C, \; K))$.

From $5.$, it immediately follows that $\delta_0 < \beta(C,K)$, since, 
as we shall see in the sequel, the bound in $5.(i)$ is smaller than 
$\beta(C,K)$. Note also that $\alpha(C,K)\geq 0$. Indeed, if 
$\alpha(C,K)<0$, then $C(C-2K)<\sqrt{C^2(C-2K)^2}$, which contradicts the 
Index theorem, since $C(C-2K)>0$.   

Observe that $\alpha(C,K)<1$ if and only if
\begin{equation}
C(C-2K)- 8 < \sqrt{C^2(C-2K)^2} \label{eq:1.16}
\end{equation}To simplify the notation, we put $t=C(C-2K)$ so that 
$(\ref{eq:1.16})$ becomes
\begin{equation}
t-8 < \sqrt{t^2 - 4 H(C,K)}. \label{eq:1.17}
\end{equation}Two cases can occur.

If $t-8<0$, there is nothing to prove since the right side member of 
($\ref{eq:1.17}$) is always positive. 

Note, before proceeding to consider the other case, that in this situation 
we want that $ \beta(C,K)>1$ in order to have at least an effective 
positive integral value for the number of nodes; but $\beta(C,K)>1$ if 
and only if $(0<) 8-t < \sqrt{t^2 - 4 H(C,K)}$. By squaring both members of 
the previous inequality, we get $4H(C,K) < 16t - 64$, which is our 
hypothesis 4.; so the upper-bound for $\delta$ is surely greater than 1. 
Moreover, the expression for such bound is the one in $5.(ii)$ and it can 
not be written in a better non-trivial form.

On the other hand, if $t-8\geq 0$, by squaring both members of 
$(\ref{eq:1.17})$, we get 
$H(C,K) < 4(C(C-2K)-4)$, which is our hypothesis 4.. Therefore, $\alpha(C,K)
 <1$; moreover, the condition $\beta(C,K) >1$ is trivially satisfied, since 
it is equivalent to $t-8 >-\sqrt{C^2(C-2K)^2}$. From $(\ref{eq:1.16})$, we 
can write

\vspace{-0.5cm}
$$   \frac{C(C-2K)+C(C-2K)-8}{8} < \frac{C(C-2K)+\sqrt{C^2(C-2K)^2}}{8},$$ 
so we can replace the bound $\delta < \beta(C,K)$ with the more 
''accessible'' one $\delta \leq \frac{C(C-2K)}{4}-1$, which is the bound in
 $5. (i)$. 

Observe that

\vspace{-0.5cm}
$$ \frac{C(C-2K)+C(C-2K)-8}{8} < \frac{C(C-2K)+\sqrt{C^2(C-2K)^2}}{8} $$ 
$$\leq \frac{C(C-2K)}{4},$$so, it is not correct to directly write $\delta<
 \frac{C(C-2K)}{4}$. Therefore, $5. (i)$ is the right approximation.

In conclusion, our numerical hypotheses contradict ($\ref{eq:1.15}$), 
therefore the assumption \par
$h^1(\Ii_N(C))$ $\neq 0$ leads to a contradiction.
\end{proof} 

\begin{remark}

\normalfont{

1) The previous theorem gives purely numerical conditions to 
deduce some informations about Severi varieties of smooth projective 
surfaces. 
In the next section, we shall discuss some interesting examples of 
projective 
surfaces 
for which the theorem can be easily applied. More precisely, we will 
consider
 smooth surfaces in $\Pt$ which are general elements of a component of the 
Noether-Lefschetz locus; for example, surfaces of general type, of 
degree $d \geq 5$, which contain a line.

Our result obviously generalizes the one of Chiantini and Sernesi. 
In their case, since $C \equiv_{num} pK_S$, $p \in \Q$ and $p \geq 2$,
we always have $ \alpha(C, K_S)=0$ and $\beta(C,K_S) = \frac{p(p-2)}{4}
K_S^2$;
this depends on the fact that $H(pK_S, K_S)=0$, for every $p$. With the 
further hypotheses that $p \in \ZZ^+$, $p$ odd, and that the 
Neron-Severi group of $S$ is $NS(S) \cong \ZZ[K_S]$, they proved 
that one can take $\delta = \frac{(p-1)^2}{4}K_S^2$. These bounds are sharp,
 at least for the general quintic surface in $\Pt$.

Furthermore, as recalled in Remark 1.2 in \cite{C-S}, in the case of rational or 
ruled surfaces (for which $CK_S<0$) or $K3$ surfaces (for which $CK_S=0$) 
if $ \mid C \mid$ is base-point-free the argument for $S=\Pd$ can be 
repeated without changes, since the line bundle $N_{\varphi}$, on 
the normalization of the nodal curve $C' \in \mid C \mid$, is non-special. 
Our result focuses on cases in which $CK_S \geq 0$ (see hypothesis $3.$), 
where the previous approach fails. 

2) One can immediately deduce that when the Hodge number is zero, i.e. 
when we are considering a divisor pair such that $(CK)^2 = C^2K^2$, 
then in the previous proof we find $\alpha(C,K)=0$ and $ \beta(C,K)= 
\frac{C(C-2K)}{4}$.  

3)Theorem $\ref{main}$ generalizes, in the case of nodes, the result in 
\cite{G-L-S II}. This will be clear after having considered the following 
examples.
}
\end{remark}

{\bf Examples:}

1) Let $S \subset \Pt$ be a general smooth quartic. We have, 
$\Oc_S(K_S) \cong \Oc_S$. Let $H$ denote the plane 
section and $D$ be the generic element of $\mid 2H \mid$. From Bertini's 
theorem, $D$ is smooth and irreducible. If $\pi: \tilde{S} \to S$ denotes 
the blow-up of $S$ in a point $P \in S$ and $E$ the associated 
exceptional divisor, then $ K_{\tilde{S}} \equiv E$, i.e. the canonical 
divisor of the blown-up surface is linearly equivalent to the exceptional 
divisor. Thus, $C \equiv 2 \pi^*(H)$ can not be ample, since $CK_{\tilde{S}}
 =0$; so, the first hypothesis in \cite{G-L-S II} does not hold. 

Nevertheless, observe that the generic element of $\mid C \mid$ is 
smooth and irreducible. Moreover, $ C- 2K_{\tilde{S}} \equiv 2 \pi^*(H) - 
2E$ so that $(C- 2K_{\tilde{S}})^2=12$, $ C(C- 2K_{\tilde{S}}) =16$, $ 
CK_{\tilde{S}}=0$, $ K_{\tilde{S}}^2=-1$, $H(C,K_{\tilde{S}})= 16$ and 
$4(C(C- 2K_{\tilde{S}})-4)=48$. Since we are in the situation $5.(i)$, 
we get $\delta \leq \frac{16}{4}-1=3$, i.e. on $\tilde{S}$, if 
$V_{\mid 2 \pi^*(H) \mid, \delta}\neq \emptyset $ and if $\delta \leq 3$, 
then it is everywhere smooth of the expected dimension.

2) Let $S$ be a smooth quintic surface in $\Pt$ which contains a line $L$. 
Denote by $\Gamma \subset S$ a plane quartic which is coplanar to $L$, so 
that $ \Gamma \equiv H-L$ ($H$ denotes the plane section). Thus,

$$H^2=5, \; HL=1, \; L^2=-3, \; H \Gamma=4, \; \Gamma^2=0 \; and \; 
\Gamma L=4.$$Choose $C \equiv 3H+L$, so that $\mid C \mid$ contains the 
curves which are residue to $\Gamma$ in the complete intersection of $S$ 
with the smooth quartic surfaces of $\Pt$ containing $\Gamma$. $\mid 3H+L 
\mid$ is base-point-free and not composed with a pencil, since $(3H+L)L=0$ 
and $3H$ is an ample divisor. By Bertini's theorems, its general member is 
smooth and irreducible; but $C$ and $C-K_S$ can not be both either ample or, 
even, nef divisors. In fact, $CL=0$ and $(C-K_S)L=(2H+L)L=-1$. Therefore, 
the result in \cite{G-L-S II} can not be applied. 

Neverthless, $CK_S= C(C-2K_S)=H(C, K_S)= 16$, $(C-2K_S)^2=3$, $K_S^2=5$, 
$4(C(C-2K_S)-4)=48$ and, since $C(C-2K_S)>8$, $\delta \leq \frac{16}{4}-1=3$.
 Thus, if $\mid 3H+L \mid $ contains some nodal, irreducible curves, then, 
if $\delta \leq 3$, $V_{\mid 3H+L \mid, \delta}$ is everywhere smooth of the 
expected dimension. 

\section*{Some results on surfaces in $\Pt$ which contain a line}

We now consider a class of examples to which our result can be easily 
applied. We shall focus on surfaces of $\Pt$ containing a line. Such 
approach can be generalized to surfaces belonging to other components 
of the Noether-Lefschetz locus.

Firstly, let $S\subset \Pt$ be a smooth quintic and $L \subset S$ a
 line. Since $p_a(L)=p_g(L)=0$, by the adjunction formula and by the fact 
that $K_S \equiv H$ we get $L^2=-3$. As before, $$ K_S^2=5, \; LH=1, \; 
L^2=-3.$$ We are interested on some results of regularity for Severi 
varieties of curves on $S$, which are residue to the line $L$ in the 
complete intersection of $S$ with a general surface of degree $a$ passing 
through the line. Thus $C \equiv aH-L$ on $S$. By straightforward 
computations, we get $$ deg(C)= (aH-L)H= 5a-1,$$ $$p_a(C)= 
\frac{5a^2+3a}{2}-1.$$ We want to find conditions on $a$ in order 
to apply our result. 

(i) $\mid C \mid$ has a smooth and irreducible general member:
For the smoothness, we have to prove that $\mid aH-L \mid$ is 
base-point-free and not composed with a pencil. Since $a \geq 1$, $aH-L = 
(a-1)H+H-L$. If $a \geq 2$, the linear system $\mid (a-1)H \mid$ can not 
have fixed intersection points on $L$. We can restrict ourselves to treat 
the behaviour of $\mid H-L \mid$ on $L$. If $\mid H-L \mid$ admitted fixed 
points on $L$, each of those points should be a tangent point for $S$ and 
the general plane of $\Pt$ passing through the line. This would imply that 
$S$ is a singular surface in such points, which contradict the hypothesis. 
Moreover, $\mid H-L \mid$ can not be composed with a pencil, since $\mid 
H-L \mid + L \subset \mid H \mid$.

For the irreducibility, we can use the fact that $C$ and $L$ are linked in 
$\Pt$ (see \cite{M}). In fact, this implies that $C$ is $projectively \;normal$ 
in $\Pt$, i.e. if we consider the exact sequence

$$ 0 \to \Ii_{C/\Pt}(\rho) \to \Oc_{\Pt}(\rho) \to \Oc_C(\rho) \to 0,$$
then $H^1(\Ii_{C/\Pt}(\rho))=0$, for each $\rho \in \ZZ$. By choosing 
$\rho=0$, we get $H^0(\Ii_{C/\Pt})=H^1(\Ii_{C/\Pt})=0$ so $H^0(\Oc_C)
\cong H^0(\Oc_{\Pt})$. This proves that $C$ is a connected curve; since we 
have already proven its smoothness, then the general member is also 
irreducible.

(ii) Numerical hypotheses: By simple computations, one 
observes that all the numerical conditions in Theorem $\ref{main}$ 
simultaneously hold if $a \geq 4$.

We can completely generalize the previous procedure to the case of a 
smooth surface of degree $d \geq 6$ which contains a line $L$. For the 
detailed computations, the reader is referred to \cite{F}.

Let $S\subset \Pt$ be such a surface and $C \equiv aH-L$, so that $$ 
deg(C)=ad-1,$$ $$ p_a(C)= \frac{ad(a+d) - 2a - d(4a+1)+3}{2}.$$ Moreover, 
$L^2= 2-d$, since $K_S \equiv (d-4)H$ and $LH=1$.

For the smoothness and the irreducibility of the general member of 
$\mid C \mid$ one can use the previous argument. Now, $K_S^2= (d-4)^2 d 
\geq 24$, because $d \geq 6$. It is not difficult to compute (see \cite{F}
 for details) that, for $6 \leq d \leq 7$, all the 
numerical hypotheses in Theorem $\ref{main}$ simultaneously hold if 
$ a \geq 2d-6$ 
(note that, for $d=5$ we obtained $a \geq 4$ so that $d=5, \; 6, \; 7$ 
behave in the same way). On the other hand, for $d \geq 8$, the condition on 
the Hodge number (i.e. hypothesis $4.$ in Theorem $\ref{main}$) determines a bound on 
$a$ which is bigger than the one determined by the other conditions, i.e. 
$2d-6$. Indeed, condition $4.$ holds if and only if $$ 4a^2d - 8a(d^2-4d+1) 
- (d^4 - 10 d^3 + 33 d^2 - 44 d + 56)>0.$$By solving this inequality, we 
find 
\begin{equation}
a >d-4+ \frac{1}{d} + \frac{1}{2}\sqrt{d^3 - 6 d^2 + d + 28 + \frac{24}{d} 
+ \frac{4}{d^2} }. \label{eq:A}
\end{equation}It is a straightforword computation to find that the right 
side member of $(\ref{eq:A})$ is bigger than $2d-6$ when $d \geq 8$. 
Therefore, in this case, all the numerical conditions of Theorem 
$\ref{main}$  
simultaneously hold if $(\ref{eq:A})$ holds.

In order to find a better expression for such a lower-bound on $a$, 
we observe that 

\vspace{-10pt}
$$ \sqrt{d^3 - 6 d^2 + d + 28 + \frac{24}{d} + \frac{4}{d^2} } <
\sqrt{d^3 - 6 d^2 + d + 32 },$$since $d \geq 8$.

We are looking for a real number $b$ such that $\sqrt{d^3 - 6 d^2 + d + 
32 } \leq \sqrt{(d \sqrt{d}- b)^2 }$. For such a value, we have 
\begin{equation}
2b\sqrt{d} \leq 6d - 1 + \frac{b^2 - 32}{d}. \label{eq:B}
\end{equation} Moreover, $(\ref{eq:A})$ becomes 
\begin{equation} 
a \geq d - 3 + \frac{d}{2} \sqrt{d} - \frac{b}{2}. \label{eq:C}
\end{equation}

Obviously, the right side member of $(\ref{eq:C})$ must be greater than 
$2d-6$ for $ d \geq 8$. Observe that this happens if and only if
\begin{equation}
d\sqrt{d} > 2d + b-6. \label{eq:D}
\end{equation} Therefore, putting $ \varphi(d):= d\sqrt{d} -2d +6$, 
from $(\ref{eq:D})$ we have $b < \varphi(d)$. The function $\varphi(d)$ is 
monotone increasing for $d \geq 2$ so, to find a uniform bound on $b$ for 
all the cases in $ d \geq 8$, it is sufficient to consider $ b < \varphi(8)$,
 i.e. $b \leq 12$. By taking into account $(\ref{eq:B})$, we find that in 
all cases a good choice is $b=9$. Thus, $(\ref{eq:A})$ can be replaced by 
$a \geq d-3 + \frac{d \sqrt{d}-9}{2}$.

Analogous computations show that, when $d \geq 5$, only condition 
2.(i) can occur, i.e. $C(C-2K_S) \geq 8$. Therefore the expression for 
the bound on the number of nodes is the one in 5.(i), which is $$ \delta 
\leq \frac{a^2d - 2a(d^2-4d+1)+ 2d -15}{4}.$$ Now, by summarizing all we 
have observed up to now, we are able to state the following
\begin{proposition}
Let $S$ be a smooth surface in $\Pt$ of degree $d \geq 5$, which contains 
a line $L$. Consider on $S$ the linear series $ \mid aH-L \mid$, with 
\begin{enumerate}
\item $a \geq 2d - 6$, if $5 \leq d \leq 7$;
\item $a \geq \lceil d-3 + \frac{d\sqrt{d}-9}{2} \rceil$, if $d \geq 8$.
\end{enumerate}
(We denote by $\lceil x \rceil$ the {\bf round-up} of the real number 
$x$, i.e. the smallest integer which is bigger than or equal to $x$). 
Suppose, also,  that for a given integer $ \delta$ the Severi variety 
$ V_{\mid aH-L \mid, \delta}$ is non-empty. Then, if $$ \delta \leq 
\frac{a^2d - 2a(d^2-4d+1)+ 2d -15}{4}, $$ the Severi variety is 
everywhere smooth of the expected dimension.
\end{proposition}

\begin{remark}
\normalfont{We want to point out that our results, in a 
certain sense, 
agree with what is stated in \cite{C-L}. In fact the authors proved the following result.} 
\end{remark}

{\bf Theorem}
{\it  Let $D$ be a reduced curve in $\Pt$ and $s$, $d$ be 
two integers such that $d \geq s+4$. Moreover, suppose that \par
i)there exists a surface $Y \subset \Pt$ of degree $s$ which contains 
$D$;\par
ii) the general element of the linear system $\mid \Oc_Y(dH - D) \mid$ 
is smooth and irreducible.\par
Denote by $S$ a general surface of $\Pt$, of degree $d$, containing $D$. 
Thus, $S$ does not contain reduced, irreducible curves $C \neq D$ of 
geometric genus $g < 1+ deg (C) \frac{(d-s-5)}{2}$. In particular, if 
$d\geq s+6$ and $p_g(D) \geq 2$, $S$ is algebraically hyperbolic.}

In the case of our proposition, $S$ is a surface of degree $d \geq 5$ 
and $D=L$, such that $L^2 = 2-d$. Thus, we can consider $s=1$, i.e. 
$Y$ is a plane containing the line $L$ and $\mid \Oc_Y(dH-L)\mid = 
\mid \Oc_{\Pd}(d-1) \mid$ which has a smooth and irreducible general 
element. Therefore, if there exists a curve $C$ of a given degree, then

$$p_g(C) \geq 2 + \frac{(d-6)}{2} deg(C) =2 + \frac{(d-6)}{2} CH.$$ 
If, moreover, $C$ is a nodal curve, then 

$$ \delta = p_a(C) - p_g(C) \leq \frac{C^2 + CK_S}{2} + 1 - 2 - 
\frac{(d-6)}{2} CH = $$ $$\frac{C^2}{2} + \frac{(d-4)}{2}CH - 
\frac{(d-6)}{2} CH  - 1= \frac{C^2}{2}+CH-1.$$ On the other hand, since 
in such cases, when all our hypotheses are satisfied, $C(C-2K_S) \geq 8$,
 $5 (i)$ determines

$$ \delta \leq \frac{C(C-2dH+8H)}{4}-1= \frac{C^2}{4}-\frac{(d-4)}{2}CH-1.$$
 Observe that $\frac{C^2}{4}-\frac{(d-4)}{2}CH-1 \leq \frac{C^2}{2}+CH-1$ 
if and only if $\frac{C^2}{4} + CH (\frac{d}{2}-1) \geq 0$. Since $d \geq 
5$ and since $C$ is big and nef (consequence of $ 1.$ and $3.$), 
this latter inequality is always strictly verified. This means that our 
bounds on $\delta$ are in the range of values, for the number of nodes, 
that are necessary for the existence of such a 
curve.

\end{document}